\documentclass[final]{siamart190516}

\usepackage{amsfonts}
\usepackage{graphicx}
\usepackage{epstopdf}
\usepackage{algorithmic}
\ifpdf
  \DeclareGraphicsExtensions{.eps,.pdf,.png,.jpg}
\else
  \DeclareGraphicsExtensions{.eps}
\fi

\newsiamremark{remark}{Remark}
\newsiamremark{hypothesis}{Hypothesis}
\crefname{hypothesis}{Hypothesis}{Hypotheses}
\newsiamthm{claim}{Claim}
\newsiamthm{assumption}{Assumption}
\newsiamremark{example}{Example}

\headers{Discrete time risk sensitive}{{\L}. Stettner}

\title{Discrete time risk sensitive control problem}

\author{{\L}ukasz Stettner\thanks{Institute of Mathematics, Polish Academy of Sciences, Warsaw, Poland,
  (\email{stettner@impan.pl});  research supported by NCN grant 2020/37/B/ST1/00463.}}

\usepackage{amsopn}

\usepackage{enumitem} 

\newcommand{\vep}{\varepsilon}

\def\bN{\mathbb{N}}

\def\ee{\mathbb{E}}

\def\prob{\mathbb{P}}
\def\tao{\tau_{\cal{O}}}

\def\o{\cal{O}}
\def\taom{\tau_{{\cal{O}}_m}}

\def\ep{\epsilon}
\def\vep{\varepsilon}
\newcommand\ind[1]{1_{#1}}

\makeatletter
\def\namedlabel#1#2{\begingroup
    #2%
    \def\@currentlabel{#2}%
    \phantomsection\label{#1}\endgroup
}
\makeatother

\begin{document}

\maketitle

\begin{abstract}
In the paper adapting Krein Rutman theory we show the existence of solutions to the long run risk sensitive control problem for controlled discrete time Markov processes over locally compact separable metric spaces.
\end{abstract}



\section{Introduction}\label{S:introduction}
Let $(X_n)$ be a discrete Markov process on $(\Omega, F,(F_n),P)$ taking values in a locally compact separable metric space $E$ endowed with Borel $\sigma$ field ${\cal E}$ and metric $\rho$. The process has a controlled transition operator (kernel) $\prob^a(x,\cdot)$, where $x\in E$ and $a\in U$ is a fixed compact set of control parameters. We shall assume that the mapping $E\times U\ni (x,a)\mapsto \prob^a(x,D)$ is Borel measurable for $D\in {\cal E}$.  Markov process is controlled using sequence $V=(a_0,a_1,\ldots, a_n\ldots)$ such that $a_n\in U$ is $F_n$ adapted. We denote by $\ee_{x}^V$ expected value corresponding to controlled Markov process with the use of sequence $V$ and starting at time $0$ from $x$. We want to maximize the following risk sensitive long run functional
\begin{equation} \label{fun1}
J_{x}(V)=\liminf_{m\to \infty}{1 \over m}  \ln\left( \ee_x^V\left[\exp\left\{\sum_{i=0}^{m-1}c(X_{i},a_{i})\right\}\right]\right),
\end{equation}
where $c:E\times U\mapsto R$ is a continuous bounded function. We shall use the following notation $\bar{c}:=\sup_{x\in E} \sup_{a\in U} c(x,a)$ and $\underline{c}:=\inf_{x\in E} \inf_{a\in U} c(x,a)$. Without loss of generality we assume that $\bar{c}\geq 0$.  We shall also assume that $c$ is not a constant function since then control problem is trivial (any strategy is optimal).
This optimization problem leads to the following Bellman equation - we are looking for a continuous function $w$ and a constant $\lambda$ such that
\begin{equation}\label{eqBell}
e^{w(x)}=\sup_{a\in U} \left[e^{c(x,a)-\lambda}\int_E e^{w(y)}\prob^a(x,dy)\right].
\end{equation}
Since the state space $E$ is locally compact separable it is also $\sigma$ compact (see section 3.5 of \cite{Conway} and Proposition 2.14 of \cite{Dudley}) and therefore  $E=\cup_{n=1}^\infty B_n$, where $B_n:=\left\{x\in E: \rho(b,x)\leq n\right\}$ with $n\in \bN$ - the set of positive integers and $b$ is a fixed point of $E$ and each ball $B_n$ is compact. Let $\mu$ be a probability measure with support $D$ contained in $B_1$.
We shall approximate solutions of \eqref{eqBell} considering controlled Markov processes on the compact sets $B_n$. Namely, let for $x\in E$, $a\in U$ and $A\in {\cal E}$
\begin{eqnarray}\label{trann}
\tilde{\prob}_n^a(x,A)&=&(1-\rho(x,B_n))^+ \prob^a(x,A\cap B_n) + \mu(A)[\prob^a(x,B_n^c) +  \nonumber \\
&& (\rho(x,B_n)\wedge 1)(1-\prob^a(x,B_n^c))],
\end{eqnarray}
where $(\cdot)^+$ stands for a positive part, while $\rho(x,B_n)\wedge 1=\min\left\{\rho(x,B_n),1\right\}$.
Whenever the controlled process leaves the set $B_{n+1}$ it starts afresh immediately from the set $D$ with initial measure $\mu$ i.e. $\tilde{\prob}_n^a(x,A)=\mu(A)$ for $x\notin B_{n+1}$. In other words,  when the process starts from $B_{n+1}^c=E\setminus B_{n+1}$ then in the next step it starts from the set $D$ with measure $\mu$. Transition probability $\tilde{\prob}_n^a$ for $x\in B_n$ is of the form  $\tilde{\prob}_n^a(x,A)=\prob^a(x,A\cap B_n) + \mu(A)\prob^a(x,B_n^c)$. In the set $B_{n+1}\setminus B_n$ we have a linear combination of measures $\mu$ and $\tilde{\prob}_n^a$ with coefficients depending on the distance from the set $B_n$, that is for $x\in B_{n+1}\setminus B_n$ we have
$\tilde{\prob}_n^a(x,A)=(1-\rho(x,B_n))^+ \left(\prob^a(x,A\cap B_n) + \mu(A)\prob^a(x,B_n^c)\right)+\rho(x,B_n)\mu(A)$.

Define inductively subtransition kernels
\begin{eqnarray}
&&\hat{\prob}_n^a(x,A):=\prob^a(x,A\cap B_n), \nonumber \\
&&(\hat{\prob}_n^{a_1,\ldots,a_m})^m(x,A):=\int_{B_n}(\hat{\prob}_n^{a_2,\ldots,a_m})^{m-1}
(y,A) {\prob}^{a_1}(x,dy)
\end{eqnarray}
for $x\in E$, $A \in {\cal E}$, $n, m\in \bN$. The kernels $(\hat{\prob}_n^{a_1,\ldots,a_m})^m(x,A)$ presents a probability of entering to the set $A$ in $m$ steps, starting from $x$, using controls $a_1,a_2,\ldots,a_m $, assuming that we do not leave the set $B_{n}$ in the meantime.

We shall impose the following assumptions
\begin{enumerate}
\item[(\namedlabel{A1}{A.1})] the mapping $E\times U\ni (x,a) \mapsto \prob^a(x,\cdot)$ is continuous in variation norm,
\item[(\namedlabel{A2}{A.2})] for each $n\in \bN$, $x\in B_n$ and open set $\o$ there is $m\in \bN$ and $a_1,\ldots,a_m\in U$ such that $(\hat{\prob}_n^{a_1,\ldots,a_m})^m(x,{\o})> 0$.
\end{enumerate}
Assumption \eqref{A1} means that whenever $U\ni a_n\to a$ and $E\ni x_n\to x$ then
\begin{equation}
\sup_{A\in {\cal E}} |\prob^{a_n}(x_n,A)-\prob^a(x,A)|\to 0
\end{equation}
as $n \to \infty$. Assumption \eqref{A2} says that for any $x\in B_n$ and open set $\o$ there is a deterministic control such that after several iterations not leaving the set $B_n$ we shall enter the set $\o$.
Let
\begin{equation}
c_n(x,a)=(1-\rho(x,B_n))^+ c(x,a) + (\rho(x,B_n)\wedge 1)\underline{c}.
\end{equation}
Clearly $c_n$ is a continuous bounded function and its value at the boundary of $B_{n+1}$ is equal to $\underline{c}$. Moreover $c_n(x,a)\uparrow c(x,a)$ and $n\to \infty$ for $x\in E$ and $a\in U$.
We have
\begin{theorem}\label{thm1} Under assumptions \eqref{A1} and \eqref{A2} for each $n\in \bN$ there exists $w_n\in C(E)$ - the space of continuous bounded functions on $E$ and a constant $\lambda_n$ such that

\begin{equation}\label{eqBelln}
e^{w_n(x)}=\sup_{a\in U} \left[e^{c_n(x,a)-\lambda_n}\int_{E} e^{w_n(y)}\tilde{\prob}_n^a(x,dy)\right]
\end{equation}
is satisfied for $x\in E$ and $\sup_{x\in B_{n+1}}w_n(x)=0$.
\end{theorem}
We want to construct a solution to \eqref{eqBell} letting in \eqref{eqBelln} $n\to \infty$ along a suitable subsequence.
Denote by $\hat{c}(x):=\sup_{a\in U}c(x,a)$.
We shall need an additional assumption
\begin{enumerate}
\item[(\namedlabel{A3}{A.3})]  for $\lambda:=\limsup_{n \to \infty} \lambda_n$  there is $\vep>0$ such that
the set \[\left\{x: \hat{c}(x)\geq \lambda-\vep\right\}\] is compact.
\end{enumerate}
\begin{remark}\label{rem1}
Assumption \eqref{A3} is satisfied when we have that
\begin{equation}
\limsup_{\rho(b,x) \to \infty} \hat{c}(x)=\tilde{c}<\lambda.
\end{equation}
In fact, then there exist $\vep>0$ and $M>0$ such that $\hat{c}(x)<\lambda-\vep$ when $\rho(b,x)\geq M$ and consequently $\left\{x: \hat{c}(x)\geq \lambda-\vep\right\} \subset \left\{x: \rho(b,x)\leq M\right\}$, which is a compact set.
\end{remark}
Alternative assumption is of the form
\begin{enumerate}
\item[(\namedlabel{A4}{A.4})]  there is $x\in E$ and $m\in \bN$, $N\in \bN$, $\vep>0$ and $a_1,a_2,\ldots,a_m \in U$ such that $\inf_{n\geq N} (\hat{\prob}_n^{a_1,\ldots,a_m})^m(x,{\cal O}_n(\vep))>0$
 with $ {\cal O}_n(\vep)=\left\{y\in E: e^{w_n(y)}>\vep\right\}$.
    \end{enumerate}
\begin{remark}
We have to avoid the case when $e^{w_n(x)}$ approaches $0$ as $n\to \infty$, as is shown in the example \ref{ex1}. Under \eqref{A4} we obtain that there is $x\in E$ such that $e^{w(x)}>0$, which by continuity and \eqref{A2} shall imply the positivity of $e^{w(y)}$ for all $y\in E$.
\end{remark}
\begin{remark}\label{rem3}
Assumption \eqref{A4} is in particular satisfied when there is a constant $L$ such that
\begin{equation}\label{addass}
\sup_{a\in U} \sup_{x,x'\in E} \sup_{B\in {\cal E}} {\prob^a(x,B)\over \prob^a(x',B)}\leq L<\infty
\end{equation}
assuming that ${0\over 0}=0$ and
\begin{equation}\label{addass'}
\sup_{B\in {\cal E}} \mu(B)\leq L \inf_{a\in U}\prob^a(b,B).
\end{equation}
In fact, we then have
\begin{equation}\label{addassp}
\sup_{a\in U} \sup_{x\in B_{n+1},x'\in B_n} \sup_{B\in {\cal E}} {\tilde{\prob}^a_n(x,B)\over \tilde{\prob}^a_n(x',B)}\leq L
\end{equation}
and consequently for $x\in B_{n+1}$ and $x'\in B_n$ we obtain

\begin{eqnarray}
&& e^{w_n(x)-w_n(x')}\leq \sup_{a\in U} e^{c(x,a)-c(x',a)}{\int_E e^{w_n(y)}\tilde{\prob}^a_n(x,dy) \over \int_E e^{w_n(y)}\tilde{\prob}^a_n(x',dy)} \leq \nonumber \\
&& \sup_{a\in U} e^{c(x,a)-c(x',a)}\left[L+{\int_E e^{w_n(y)}(\tilde{\prob}^a_n(x,dy)-L \tilde{\prob}^a_n(x',dy)) \over \int_E e^{w_n(y)}\tilde{\prob}^a_n(x',dy)}\right]\leq \nonumber \\
&& \sup_{a\in U} e^{c(x,a)-c(x',a)}L.
\end{eqnarray}
Then $\sup_{x\in B_{n+1}}w_n(x)-\inf_{x'\in B_n}w_n(x') = 0-\inf_{x'\in B_n}w_n(x')<K$ with $K$ chosen independently on $n\in \bN$. Therefore $\inf_{x'\in B_n}w_n(x')>-K$ and $B_n\subset {\cal O}_n(e^{-K})$. It is clear that there is $N\in \bN$, $x\in E$ and $a\in U$ such that $\prob^a(x,B_N)>0$. Then $\inf_{n\geq N} \hat{\prob}^a_n(x,{\cal O}_n(e^{-K}))>0$ and assumption \eqref{A4} is satisfied.
\end{remark}
We have
\begin{theorem}\label{thm2} Under \eqref{A1}, \eqref{A2}, when $E$ is not compact and either \eqref{A3} or \eqref{A4} is satisfied, there is a continuous function $w$ such that $\sup_{x\in E}w(x)=0$ and the pair $(w,\lambda)$ forms a solution to the equation \eqref{eqBell}.
\end{theorem}
In the case when $E$ is compact the solution to \eqref{eqBell} can be obtained in a simpler way, see Proposition  \ref{propex1}.
The problem of existence of solutions to the Bellman equation \eqref{eqBell} in locally compact separable state spaces was studied first in \cite{DiMasi 1999} under uniform ergodicity assumption and equivalence of transition probabilities.  Then in \cite{DiMasi 2000} uniform ergodicity was assumed together with under small risk factor (sufficiently small span norm of $c$). In the paper \cite{DiMasi 2007} using splitting arguments uniform ergodicity was replaced by minorization property together with sufficiently small norm of $c$. Control problem with cost functional \eqref{fun1} in general state space was also studied in \cite{Jas2007}, where instead of equation a Bellman inequality was considered and the approach was based on certain assumption on discounted approximation.
The use of Krein Rutman theorem (see \cite{KreinRut} or \cite{Bonsall}, \cite{Arap})  to risk sensitive problems appeared first in the paper \cite{AnanB}, where existence of solutions to the risk sensitive Bellman equation over compact state space was considered. Similar ideology was the used in the paper \cite{Bis}, where countable state space was studied and an approximation solution to the Bellman equation over countable state space was considered, using sequence solutions in finite state spaces.  There is a long list of contributors studying the problem with finite or countable state spaces (see \cite{Cadena}, \cite{CadenaH}, \cite{Bis} and references therein). Continuous time log run risk sensitive problems also were studied intensively (see e.g. \cite{ArapBis}, \cite{ArapBisB}, \cite{Guo}, or \cite{PitSte2019}).
In the paper we consider discrete time long run risk sensitive control problem over locally compact separable metric space. We use an approximation of the solution to the risk sensitive Bellman equation based on solutions to Bellman equations restricted to an increasing sequence of balls. We obtain solutions to risk sensitive Bellman equation in balls using famous Krein Rutman theorem (see \cite{KreinRut} or \cite{Bonsall}, \cite{Arap}). The novelty of the paper comparing to \cite{AnanB} is that we consider general locally compact separable state space and use more general assumptions. We also consider an example to point out possible problems with above mentioned approximation.   We study first the case with maximization of the risk sensitive functional and then sketch the differences in the case of minimization. The problem with maximization corresponds to the studies of asymptotic behaviour of power utility function, while minimization is closely related to risk sensitive maximization (see e.g. \cite{PitSte2022}).

\section{Solution to the equation \eqref{eqBelln}}

Let for $f \in C(E)$
\begin{equation}
T_nf(x):= \sup_{a\in U}\left[e^{c_n(x,a)}\int_E f(y)\tilde{\prob}_n^a(x,dy)\right].
\end{equation}
For $x\in B_{n+1}^c$ we have that $T_nf(x)=e^{{\underline{c}}}\mu(f)$.
We have
\begin{proposition}\label{prop1}
Under \eqref{A1} operator $T_n$ is completely continuous in the space $C(E)$ i.e. it transforms bounded subsets of $C(E)$ into compact subsets of $C(E)$. For function $f\in C(E)$ and $x,x'\in E$ we have
\begin{eqnarray}\label{equicont}
|T_nf(x)-T_nf(x')|&\leq & \sup_{a \in U}|e^{c_n(x,a)}-e^{c_n(x',a)}|\|f\|+ \nonumber \\
&& e^{\bar{c}}\|f\|\left[2\rho(x,x')+3\sup_{a\in U}\|\prob^a(x,\cdot)-\prob^a(x',\cdot)\|_{var}\right]
\end{eqnarray}
where $\|\cdot\|_{var}$ stands for the total variation norm.
\end{proposition}
\begin{proof}
Notice that
\begin{eqnarray}\label{eqn1}
&&|T_nf(x)-T_n(f(x')|\leq \sup_{a\in U}\left[|e^{c_n(x,a)}-e^{c_n(x',a)}|\|f\|+ \right. \nonumber \\
&& \left. e^{c_n(x',a)}|\tilde{\prob}_n^af(x)-
\tilde{\prob}_n^af(x')|\right].
\end{eqnarray}
Now
\begin{eqnarray}\label{eqn2}
&&|\tilde{\prob}_n^af(x)-\tilde{\prob}_n^af(x')|\leq \sup_{a\in U} |\int_{E}f(y){\prob}_n^a(x,dy)\left[(1-\rho(x,B_n))^+ - \right. \nonumber \\
&& \left. (1-\rho(x',B_n))^+\right] +(1-\rho(x',B_n))^+ \left[\int_{E} f(y) ({\prob}_n^a(x,dy)- {\prob}_n^a(x',dy))\right] + \nonumber \\
&&\mu(f)\left[{\prob}_n^a(x,B_n^c)-{\prob}_n^a(x',B_n^c)+ \left(\rho(x,B_n)\wedge 1 - \rho(x',B_n)\wedge 1\right)(1-\right. \nonumber \\
&& \left. \tilde{\prob}_n^a(x,B_n^c))+(\rho(x',B_n)\wedge 1) ({\prob}_n^a(x',B_n^c)-{\prob}_n^a(x,B_n^c))\right]|.
\end{eqnarray}
Taking into account that
\begin{equation}
|(1-\rho(x,B_n))^+ - (1-\rho(x',B_n))^+|\leq \rho(x,x')
\end{equation}
and
\begin{equation}
|\rho(x,B_n)\wedge 1 - \rho(x',B_n)\wedge 1| \leq \rho(x,x')
\end{equation}
we obtain
\begin{equation}\label{eqn3}
|\tilde{\prob}_n^af(x)-\tilde{\prob}_n^af(x')|\leq 3\|f\| \sup_{a\in U}\|\tilde{\prob}_n^a(x,\cdot)-\tilde{\prob}_n^a(x',\cdot)\|_{var} + 2 \rho(x,x')\|f\|.
\end{equation}
Now from \eqref{eqn1} and \eqref{eqn3} we obtain \eqref{equicont}. Since the value of $T_nf(x)$ for $x \in B_{n+1}^c$ is equal to $e^{\underline{c}}\mu(f)$ we have that
for uniformly bounded functions $f$ the values of $T_nf$ are equicontinuous and by Ascoli Arzela theorem (see Theorem 7.40 of \cite{Royden}) closure of the image of the operator $T_n$ of any bounded set in $C(E)$ is compact.
\end{proof}
Therefore by the Krein Rutman theorem (see Theorem \ref{thmA1}) we obtain
\begin{corollary}\label{cor1} Under \eqref{A1} there is $v_n\in C_+(E)$ - the class of continuous, bounded, nonnegative functions on $E$ with $\|v_n\|=1$ and a positive constant
 $\tilde{\lambda}_n$  such that for $x\in E$, $n\in \bN$
 \begin{equation}\label{Bellng}
 T_nv_n(x)=\tilde{\lambda}_n v_n(x).
 \end{equation}
 Furthermore $v_n(x)=(\tilde{\lambda}_n)^{-1}e^{\underline{c}}\mu(\nu_n)$ for $x\notin B_{n+1}$.
 \end{corollary}
 \begin{proof} We use Theorem \ref{thmA1} with $K=C_+(E)$ noticing that  $T_n 1(x)\geq e^{\underline{c}}1(x)$, where $1(x)=1$ for $x\in E$. Since for $x\notin B_{n+1}$ we have that $T_nf(x)=e^{\underline{c}}\mu(f)$ we obtain that $v_n$ is constant on $B_{n+1}^c$.
 \end{proof}
 Let $\lambda_n:=\ln \tilde{\lambda}_n$. We would like to have that $v_n(x)>0$ for $x\in E$. For this purpose we need assumption \eqref{A2}. We have the following result which completes the proof of Theorem \ref{thm1}
 \begin{lemma}\label{lem1}
 Under \eqref{A1} and \eqref{A2} there are $w_n\in C(E)$ and $\lambda_n$ such that
 $\sup_{y\in E}w_n(y)=\sup_{y\in B_{n+1}}w_n(y)=0$ and
 equation \eqref{eqBelln} is satisfied.
 \end{lemma}
 \begin{proof} Let $\o$ $=\left\{y\in E: v_n(y)>\alpha\right\}$. Since $\|v_n\|=1$ and $v_n$ is nonnegative this set is nonempty for $\alpha<1$. By assumption \eqref{A2} for any $x\in B_n$ there is $m\in \bN$ and $a_1,a_1,\ldots,a_m\in U$ such that $(\hat{\prob}_n^{a_1,\ldots,a_m})^m(x,{\o})>0$. Iterating \eqref{Bellng}
 we obtain
 \begin{equation}\label{vn positive}
(\tilde{\lambda}_n)^{m} v_n(x)\geq\alpha e^{m {\underline{c}}} (\hat{\prob}_n^{a_1,\ldots,a_m})^m(x,{\o})> 0.
\end{equation}
 Therefore $v_n(x)>0$ for $x\in B_n$ and  $w_n(x)=\ln v_n(x)$ is well defined for $x\in B_n$ and it is continuous and therefore bounded on compact sets of $B_n$. Moreover $v_n$ is a positive constant on $B_{n+1}^c$, since measure $\mu$ has a support contained in $B_1\subset B_n$ and $\mu(v_n)>0$. It remains to show that $w_n(x)$ is well defined for $x\in B_{n+1}\setminus B_n$. Let $G=\left\{x: v_n(x)=0\right\}$. Clearly $G \subset B_{n+1}\setminus B_n$, and by \eqref{Bellng}, taking into account nonnegativity of $v_n$  we have that $\tilde{\prob}_n^a(x,G)=1$ for any $a\in U$ and $x\in G$. On the other hand  from \eqref{trann} we have $\tilde{\prob}_n^a(x,G)=0$, a contradiction.
 Consequently $w_n\in C(E)$ together with $\lambda_n$ satisfy \eqref{eqBelln}. Function $w_n$ is constant equal to $\underline{c}+\ln \mu(v_n)- \ln \tilde{\lambda}_n$ on $B_{n+1}^c$ and by continuity $w_n(x)=\underline{c}+\ln \mu(v_n)- \ln \tilde{\lambda}_n$ also on the boundary of $B_{n+1}$, which completes the proof.
 \end{proof}

 \begin{proposition}\label{propex1}
 When $E$ is compact then under \eqref{A1} and the following version of \eqref{A2}
 \begin{enumerate}
 \item[(\namedlabel{A2*}{A.2*})] for each  $x\in E$ and open set $\o$ there is $m\in \bN$ and $a_1,\ldots,a_m\in U$ such that $({\prob}^{a_1,\ldots,a_m})^m(x,{\o})> 0$.
\end{enumerate}
there is $w\in C(E)$ and $\lambda$ which form a solution to the equation \eqref{eqBell}.
\end{proposition}
\begin{proof} We use Krein Rutman theorem (Theorem \ref{thmA1}) as in Corollary \ref{cor1} and then follow the proof of Lemma \ref{lem1} using \eqref{A2*}.
\end{proof}
 Using stochastic control interpretation of the Bellman equation \eqref{eqBelln} we easily obtain
 \begin{corollary}\label{cor2}
Under \eqref{A1} and \eqref{A2} we have the following bounds for $\lambda_n$
\begin{equation}\label{bounds}
\underline{c}\leq \lambda_n \leq \bar{c}.
\end{equation}
\end{corollary}
\begin{proof} It follows directly from suitable version of Proposition 1.1 of \cite{DiMasi 1999} that
\begin{equation}
\lambda_n=\sup_V \liminf_{m\to \infty} {1\over m} \ln \left(\tilde{\ee}_x^V\left[\exp\left\{\sum_{i=0}^{m-1}c_n(X_i,a_i)\right\}\right]\right)
\end{equation}
where $\tilde{\ee}_x^V$ is the expectation with respect to the controlled Markov process starting from $x$  with control $V$ and transition probabilities defined in \eqref{trann}. The bounds \eqref{bounds} are then immediate.
\end{proof}

\begin{remark}
Notice that upper bound for $\lambda_n$ is $\bar{c}$ and this is satisfied even without assumption \eqref{A2}. In fact iterating \eqref{Bellng} we obtain for $m\in \bN$
\begin{eqnarray}\label{eqiter}
&&v_n(x)=\sup_V \tilde{\ee}_x^V\left[\exp\left\{\sum_{i=0}^{m-1}(c_n(X_i,a_i)-\lambda_n))\right\}v_n(X_m)\right]\leq \nonumber \\
&& \sup_V \tilde{\ee}_x^V\left[\exp\left\{\sum_{i=0}^{m-1}(c_n(X_i,a_i)-\lambda_n)\right\}\right]
\end{eqnarray}
and consequently
\begin{equation}
\lambda_n\leq\sup_V \liminf_{m\to \infty} {1\over m} \ln \left(\tilde{\ee}_x^V\left[\exp\left\{\sum_{i=0}^{m-1}c_n(X_i,a_i)\right\}\right]\right)\leq \bar{c}.
\end{equation}
\end{remark}
\begin{remark}
One can easily notice that under \eqref{A1} and \eqref{A2} value $\lambda_n$ coincides with the type of an optimally controlled semigroup (see \cite{Hille} Sec. 10.2), namely
\begin{equation}
\lambda_n=\sup_V \liminf_{m\to \infty} {1\over m} \ln \sup_{x\in E} \left(\tilde{\ee}_x^V\left[\exp\left\{\sum_{i=0}^{m-1}c_n(X_i,a_i)\right\}\right]\right).
\end{equation}
Consequently $\lambda_n$ is unique. On the other hand $w_n$ as a solution to \eqref{eqBelln} is not unique. In particular any $w_n(x)+d$, where $d$ is a constant is also a solution (together with $\lambda_n$) to the equation \eqref{eqBelln}. Notice furthermore that although $c_n(x,a)\leq c_{n+1}(x,a)$ the values of $\lambda_n$ are not necessarily increasing.
\end{remark}

\section{Solution to the equation \eqref{eqBell}}
In this section we solve equation \eqref{eqBell} using suitable subsequence of $n\to \infty$ in the equation \eqref{eqBelln}.
We start with the following
\begin{lemma}\label{lem2} Under \eqref{A1} we have that
$\tilde{\prob}_n^a(x,\cdot)$ converges uniformly in $x$ from compact subsets of $E$ and $a\in U$, in the total variation norm to $\prob^a(x,\cdot)$.
\end{lemma}
\begin{proof}
By \eqref{A1} for a given compact set compact ball $B_k$ the set of measures $\left\{\prob^a(x,\cdot), x\in B_k, a\in U\right\}$ is compact in weak topology. Therefore by Theorem 5.2 of \cite{Bil}  for a given $\vep>0$ and there is $m>k$ such that $\sup_{x\in B_k} \sup_{a\in U} \prob^a(x,B_m^c)\leq \vep$.  For $x\in B_k$, $a\in U$ we have for $n\geq m$
\begin{eqnarray}
&&|\tilde{\prob}_n^a(x,A)-{\prob}^a(x,A)| \leq (1-\rho(x,B_n))^+ |\prob^a(x,A\cap B_n)-\prob^a(x,A)| + \nonumber \\ &&(1-(1-\rho(x,B_n))^+)\prob^a(x,A)+
\mu(A\cap D)[\prob^a(x,B_n^c) + \nonumber \\
&& (\rho(x,B_n)\wedge 1)(1-\prob^a(x,B_n^c))]\leq
 \prob^a(x,A\cap B_n^c)(1-\rho(x,B_n))^+  +\nonumber \\
 && (1-(1-\rho(x,B_n))^+)+
 \prob^a(x,B_n^c)+\rho(x,B_n) \leq 2 \vep.
\end{eqnarray}
Since $\vep$ could be chosen arbitrarily small we have the claim.
\end{proof}

We shall now prove Theorem \ref{thm2}. By \eqref{bounds} we can choose a  subsequence $n_k$  such that $\lambda_{n_k}\to \lambda=\limsup_{n\to \infty}\lambda_n$, as $k\to \infty$.
For $x,x'\in E$, taking into account that $\|e^{w_n}\|=1$ from \eqref{equicont} we have
\begin{eqnarray}\label{equicontn}
|e^{w_n(x)}-e^{w_n(x')}|&\leq & \sup_{a \in U}|e^{c_n(x,a)}-e^{c_n(x',a)}|+ \nonumber \\
&& e^{\bar{c}}\left[2\rho(x,x')+3\sup_{a\in U}\|\prob^a(x,\cdot)-\prob^a(x',\cdot)\|_{var}\right].
\end{eqnarray}
This means that sequence of functions $e^{w_n(x)}$ is equicontinuous and bounded. Therefore by Ascoli Arzela theorem (see \cite{Royden} Theorem 7.40) there is a further subsequence of $(n_k)$, for simplicity denoted by $(n_k)$ and function $z\in C(E)$ such that $e^{w_{n_k}(x)}\to z(x)$, as $k\to \infty$, uniformly on compact subsets of $E$.
Since $\|e^{w_n}\|=1$ we therefore have that $\|z\|\leq 1$.  For a given $\vep>0$ there is a compact set $B_m$ such that for $n\in \bN$ and fixed $x\in E$
\begin{equation}\label{o1}
\sup_{a\in U} \tilde{\prob}_n^a(x,B_m^c)\leq \vep.
\end{equation}
Consequently for $n\geq m$ we obtain
\begin{eqnarray}\label{o2}
&&\sup_{a\in U}|e^{c_n(x,a)-\lambda_n}\int_{E} e^{w_n(y)}\tilde{\prob}_n^a(x,dy) - e^{c(x,a)-\lambda}\int_{E} z(y){\prob}^a(x,dy)|   \leq \nonumber \\
&& \sup_{a\in U} |e^{c_n(x,a)-\lambda_n}-e^{c(x,a)-\lambda}|\int_{E} e^{w_n(y)}\tilde{\prob}_n^a(x,dy)+ \nonumber \\
&&+
\sup_{a\in U}e^{c(x,a)-\lambda}| \int_{E} e^{w_n(y)}\tilde{\prob}_n^a(x,dy)-\int_{E} z(y){\prob}^a(x,dy)|\leq \nonumber \\
&& \sup_{a\in U}|e^{c_n(x,a)-\lambda_n}-e^{c(x,a)-\lambda}| +
\sup_{a\in U} e^{\|c\|-\lambda}\left[\int_{B_m} |e^{w_n(y)}-z(y)|\tilde{\prob}_n^a(x,dy)+ \right. \nonumber \\
&& \left. \tilde{\prob}_n^a(x,B_m^c) +
 |\int_E z(y)(\tilde{\prob}_n^a(x,dy)-{\prob}^a(x,dy))|\right] \leq \nonumber \\
 && \sup_{a\in U}|e^{c_n(x,a)-\lambda_n}-e^{c(x,a)-\lambda}|+ \nonumber \\
 && e^{\|c\|-\lambda}\left[\sup_{y\in B_m}|e^{w_n(y)}-z(y)|+\vep+|\int_E z(y)(\tilde{\prob}_n^a(x,dy)-{\prob}^a(x,dy))|\right].
 \end{eqnarray}
 Letting $n_k\to \infty$, since $\vep$ can be arbitrarily small we obtain that
 \begin{equation}\label{o3}
 \sup_{a\in U}|e^{c_{n_k}(x,a)-\lambda_{n_k}}\int_{E} e^{w_{n_k}(y)}\tilde{\prob}_{n_k}^a(x,dy) -  e^{c(x,a)-\lambda}\int_{E} z(y){\prob}^a(x,dy)|\to 0
 \end{equation}
 and therefore
 \begin{equation}\label{o4}
 \lim_{k\to \infty} e^{w_{n_k}(x)}= \sup_{a\in U} e^{c(x,a)-\lambda}\int_{E} z(y){\prob}^a(x,dy)=z(x).
 \end{equation}
 It remains to show that $w(x)=\ln z(x)$ is well defined and $w\in C(E)$. Notice first that  $w_n\not\equiv 0$ since otherwise by \eqref{eqBelln} we have $w_n(x)=\sup_{a\in U} c_n(x,a)-\lambda_n=0$ for  $x\in E$ and letting $k\to \infty$ in $w_{n_k}(x)$ we obtain that $c\equiv \lambda$ is a constant function, which we excluded.

 Let $\hat{x}_{n_k}\in E$ be such that $w_{n_k}(\hat{x}_{n_k})=0$. Such $\hat{x}_{n_k}$ exists since by Lemma \ref{lem1} $\sup_{x\in B_{n+1}}w_n(x)=0$.
 Let $\hat{c}_n(x)=\sup_{a\in U}c_n(x,a)$. Then using \eqref{eqBelln} there is $\hat{a}_{n_k}\in U$ such that
\begin{equation}\label{o5}
1=e^{w_{n_k}(\hat{x}_{n_k})}=e^{c_{n_k}(\hat{x}_{n_k},\hat{a}_{n_k})-\lambda_{n_k}}\int_{E} e^{w_{n_k}(y)}\tilde{\prob}_{n_k}^{\hat{a}_{n_k}}(\hat{x}_{n_k},dy)\leq  e^{\hat{c}_{n_k}(\hat{x}_{n_k})-\lambda_{n_k}}
\end{equation}
and therefore $\hat{c}_{n_k}(\hat{x}_{n_k})\geq \lambda_{n_k}$.
Since $c_n(x,a)\leq c(x,a)$ for $x\in E$ and $a\in U$ we have also $\hat{c}_n(x)\leq \hat{c}(x)$ for $x\in E$, so that
$$\hat{x}_{n_k}\in \left\{x: \hat{c}_{n_k}(x)\geq \lambda_{n_k}\right\}\subset \left\{x: \hat{c}(x)\geq \lambda - (\lambda -\lambda_{n_k})\right\}.$$
Consequently for $\vep>0$ from \eqref{A3} we have a further subsequence of $(n_k)$, for simplicity denoted again as $(n_k)$ such that for  $k\geq k_0$ we have
$|\lambda-\lambda_{n_k}|\leq \vep$ and $\hat{x}_{n_k}\in \left\{x: \hat{c}(x)\geq \lambda -\vep\right\}$, which is compact.
 Since $e^{w_{n_k}(x)}\to z(x)$ uniformly in $x$ from compact sets we have that $\|z\|=1$. Using equation

\begin{equation}\label{o6}
z(x)=\sup_{a\in U} e^{c(x,a)-\lambda}\int_{E} z(y){\prob}^a(x,dy)
\end{equation}
and assumption \eqref{A2} in a similar way as in the proof of Lemma \ref{lem1} we obtain that $z(x)>0$ for $x\in E$. This completes the proof of Theorem \ref{thm2} under \eqref{A3}.

Assume now \eqref{A4}. Let for $a_1,a_2,\ldots,a_m \in U$, $m\in \bN$, $A\in {\cal E}$,
\begin{equation}
(\tilde{\prob}_n^{a_1,\ldots,a_m})^m(x,A):=\int_{E}(\tilde{\prob}_n^{a_2,\ldots,a_m})^{m-1}
(y,A) \tilde{\prob}^{a_1}(x,dy).
\end{equation}
Clearly $(\hat{{\prob}}_n^{a_1,\ldots,a_m})^m(x,A)\leq (\tilde{\prob}_n^{a_1,\ldots,a_m})^m(x,A)$ for $A\in {\cal E}$, $x\in B_n$ and $m\in \bN$, so that
subtransition kernels $(\hat{{\prob}}_n^{a_1,\ldots,a_m})^m$ are dominated by transition kernels $(\tilde{\prob}_n^{a_1,\ldots,a_m})^m$ for $x\in B_n$. Consequently as in \eqref{eqiter} we obtain
\begin{eqnarray}\label{prA3}
&& v_n(x)=\sup_V \tilde{\ee}_x^V\left[\exp\left\{\sum_{i=0}^{m-1}(c_n(X_i,a_i)-\lambda_n))\right\}
v_n(X_m)\right] \geq \nonumber \\
&& \vep e^{m\underline{c} - m \lambda_n} (\hat{\prob}_n^{a_1,\ldots,a_m})^m(x,{\cal O}_n(\vep)).
\end{eqnarray}
Letting $n_k\to \infty$, by \eqref{A4} and taking into account inequality $\underline{c}\leq \lambda_n\leq \bar{c}$ (see Corollary \ref{cor2}) we obtain that $z(x)>0$ for $x\in B_{n_k}$ and therefore also for $x\in E$. This completes the proof under \eqref{A4}.

\begin{remark}
Without assumptions \eqref{A3} or \eqref{A4} it may happen that $w_n(x)\to -\infty$ for $x\in E$ as $n\to \infty$, as is shown in the example below.
\end{remark}
\begin{example}
Let $E=\bN$. The dynamics of $(X_n)$ is deterministic $x_{n+1}=x_n+1$, $x_0\in E$ and for simplicity consider non controlled case so that $c(x,a)=c(x)=\bar{c}-{1\over x}$. We assume that $B_n:=\left\{1,2,\ldots,n+1\right\}$, $D:=\left\{1\right\}$ and $\mu=\delta_{\left\{1\right\}}$. Then for $i\in B_n$ we have
\begin{equation}\label{ex1}
e^{w_n(i)}=e^{\sum_{j=i}^{n+1}(c(j)-\lambda_n)+w_n(n+2)}
\end{equation}
and clearly $w_n(n+2)=w_n(1)$. Letting $i=1$ in \eqref{ex1} we obtain that
\begin{equation}
\lambda_n={1\over n+1} \sum_{j=1}^{n+1}c(j)=\bar{c}-{1 \over n+1}\sum_{j=1}^{n+1}{1\over  j}.
\end{equation}
Let $k=\inf\left\{j\geq 1: c(j)\geq \lambda_n\right\}$. One can calculate that
$k=\left[{n+1 \over \sum_{j=1}^{n+1} {1\over j}}\right]+1$, where $\left[\cdot\right]$ stands for the integer part. Clearly $k\to \infty$ as $n\to \infty$.

Since $e^{w_n(k)}=e^{\sum_{j=k}^{n+1} (c(j)-\lambda_n)+w_n(1)}$ we see that $w_n$ attains its maximum at $k$. We assume that $w_n(k)=0$.
Then
\begin{equation}
w_n(1)=-\sum_{j=k}^{n+1}(c(j)-\lambda_n)={k\over n+1} \sum_{j=k}^{n+1} {1\over j}- {n+1-k\over n+1}\sum_{j=1}^{k-1} {1\over j}.
\end{equation}
Since ${n+1-k \over n+1}\to 1$ as $n\to \infty$ and ${k\over n+1}\sum_{j=k}^{n+1}{1\over j}\leq 1$ we have that $w_n(1)\to -\infty$ as $n\to \infty$. Similarly one can show that $w_n(j)\to -\infty$, as $n\to \infty$ for $j\in E$.
\end{example}

\begin{remark}
Notice that the solution $w$ to \eqref{eqBell} may be unbounded from below i.e. we may have that $\liminf_{\|x\|\to \infty} w(x)=-\infty$. This is the case when outside of a compact set $B_m$ function $c(x,a)=\underline{c}$, and $\lambda>\underline{c}$. When process starts far away from the set $B_m$ at each step the negative value $\underline{c}-\lambda$ is incurred untill the process enters $B_m$. Then iterating \eqref{eqBell}
we have
\begin{eqnarray}
&&e^{w(x)}=\sup_{V} \ee_x^V\left[\exp\left\{{T_m\wedge n}(\underline{c}-\lambda)+w(X_{T_m\wedge n})\right\}\right]\leq \nonumber \\
&& \sup_{V} \ee_x^V\left[\exp\left\{{T_m\wedge n}(\underline{c}-\lambda)\right\}\right]
\end{eqnarray}
with $T_m=\inf\left\{i\geq 0: X_i\in B_m\right\}$
and letting $\|x\|\to \infty$, assuming that controlled process is not shifted immediately to $B_m$ we have that $w$ is unbounded from below.
\end{remark}
Next Proposition explains importance of the Bellman equation \eqref{eqBell}.
\begin{proposition}\label{prop2}
Assume \eqref{A1}, \eqref{A2}, and \eqref{A3} or \eqref{A4}.  We have that $\lambda\leq \inf_{x\in E} J_x(\hat{V})$. If for a given $x\in E$, any $\vep>0$ and any control $V$ there is $n\in \bN$ such that
\begin{equation}\label{exass}
\limsup_{m\to \infty}{1\over m} \ln {\ee_x^V\left(\exp\left\{\sum_{i=0}^{m-1}c(X_i,a_i)\right\}\right) \over
\ee_x^V\left(1_{B_n}(X_m)\exp\left\{\sum_{i=0}^{m-1}c(X_i,a_i)\right\}\right)}\leq \vep
\end{equation}
then $\lambda=\sup_V J_x(V)$ and an optimal control is in the form of sequence $\hat{V}=(\hat{u}(X_i))$, where $\hat{u}$ is a Borel measurable selector of the right hand side of the equation \eqref{eqBell}.
\end{proposition}
\begin{proof}
Iterating \eqref{eqBell} we obtain that (see e.g. \cite{DiMasi 1999})
\begin{equation}\label{itBell}
e^{w(x)}=\sup_V \ee_x^V\left[\exp\left\{\sum_{i=0}^{m-1}(c(X_i,a_i)-\lambda)+w(X_m)\right\}\right].
\end{equation}
Then for $\hat{V}$  taking into account that $w\leq 0$ we have
\begin{eqnarray}\label{itBellest}
&&\lambda={1\over m} \ln \left\{\ee_x^{\hat{V}}\left[\exp\left\{\sum_{i=0}^{m-1}c(X_i,\hat{u}(X_i))+w(X_m)\right\}\right]\right\}
-{w(x)\over m}\leq \nonumber \\
&&{1\over m} \ln \left\{\ee_x^{\hat{V}}\left[\exp\left\{\sum_{i=0}^{m-1}c(X_i,\hat{u}(X_i)))\right\}\right]\right\}
-{w(x)\over m}
\end{eqnarray}
and letting $m\to \infty$ we obtain that $\lambda\leq \inf_{x\in E}J_x(\hat{V})$.

Choose now $x\in E$ for which \eqref{exass} is satisfied.  Therefore for any strategy $V$ we have
\begin{eqnarray}\label{iine}
&&\lambda\geq {1\over m} \ln \left\{\ee_x^V\left[1_{B_n}(X_m)\exp\left\{\sum_{i=0}^{m-1}c(X_i,a_i)\right\}\right]
e^{\inf_{y\in B_n}w(y)}+\right. \nonumber \\
&&\left. \ee_x^V\left[1_{B_n^c}(X_m)e^{-m\|c\|+w(X_m)}\right]\right\}-{w(x)\over m}\geq \nonumber \\
&& {1\over m} \ln \left\{\ee_x^V\left[1_{B_n}(X_m)\exp\left\{\sum_{i=0}^{m-1}c(X_i,a_i)\right\}\right]\right\}+
{\inf_{y\in B_n}w(y)\over m}-{w(x)\over m}
\end{eqnarray}
Letting $m\to \infty $ by \eqref{exass} we obtain
\begin{equation}
\lambda \geq \limsup_{m\to \infty}{1\over m} \ln \ee_x^V\left(\exp\left\{\sum_{i=0}^{m-1}c(X_i,a_i)\right\}\right)-\vep
\end{equation}
Since $\vep$ can be chosen arbitrarily small and $\lambda\leq J_x(\hat{V})$ we have
\begin{equation}
\lambda = \limsup_{m\to \infty}{1\over m} \ln \ee_x^{\hat{V}}\left(\exp\left\{\sum_{i=0}^{m-1}c(X_i,a_i)\right\}\right),
\end{equation}
which completes the proof.
\end{proof}
\begin{remark}
The condition \eqref{exass} says that controlled process enters the set $B_n^c$ with probability which decreases sufficiently fast. It is satisfied when for any control $V$
there is $n\in \bN$ such that
\begin{eqnarray}
&&\limsup_{m\to \infty}{1\over m} \ln \ee_x^V\left(1_{B_n}(X_m)\exp\left\{\sum_{i=0}^{m-1}c(X_i,a_i)\right\}\right)= \nonumber \\
&&\limsup_{m\to \infty}{1\over m} \ln \ee_x^V\left(\exp\left\{\sum_{i=0}^{m-1}c(X_i,a_i)\right\}\right).
\end{eqnarray}
\end{remark}
\begin{remark}\label{remm}
One can see that under the assumptions to Proposition \ref{prop2} together with \eqref{exass} the value $\lambda$ is defined in a unique way. The function $w$ is not unique since adding any constant to $w$ we also obtain a solution to \eqref{eqBell}. On the other hand without assumption \eqref{exass} from \eqref{itBellest} we get only $\lambda\leq \inf_{x\in E}J_x(\hat{V})$ which is  an upper estimate for $\lambda$. We have defined $\lambda$ as $\limsup_{n\to \infty} \lambda_n$, which is unique. However we can think also about $\lambda$ as a limit of any subsequence of $\lambda_n$ (even smaller than $\limsup_{n\to \infty} \lambda_n$) such that \eqref{A3} is satisfied. Then we also get a solution to \eqref{eqBell}. Consequently in general case $\lambda$ as a solution to \eqref{eqBell} may not be unique. From the control point of view the optimal value of the cost functional \eqref{fun1} may depend on initial state as is shown in the Example 1 considered in \cite{Jas2007} and does not satisfy \eqref{eqBell} for any function $w$.
\end{remark}

\section{Problem with minimization}
We shall now minimize the cost functional
\begin{equation} \label{mfun1}
J_{x}^m(V)=\limsup_{m\to \infty}{1 \over m}  \ln\left( \ee_x^V\left[\exp\left\{\sum_{i=0}^{m-1}c(X_{i},a_{i})\right\}\right]\right).
\end{equation}
In the case of such multiplicative functional the problems with maximization and minimization are different. Fortunately we can adapt some technics of sections 1-3 to the minimization case but we shall need a stronger assumptions.
The corresponding Bellman equation is of the form
\begin{equation}\label{meqBell}
e^{w(x)}=\inf_{a\in U} \left[e^{c(x,a)-\lambda}\int_E e^{w(y)}\prob^a(x,dy)\right]
\end{equation}
and we are looking for a continuous function $w$ and constant $\lambda$.
To solve \eqref{meqBell} we approximate it using controlled Markov processes with transition probabilities of the form \eqref{trann}.

We consider the following stronger version of the assumption \eqref{A2}
\begin{enumerate}
\item[(\namedlabel{A2'}{A.2'})] for each $n\in \bN$, $x\in B_n$ and open set $\o$ there is $m\in \bN$  such that
    $$\inf_{V_m}(\hat{\prob}_n^{V_m})^m(x,{\o})> 0,$$
\end{enumerate}
where $V_m$ is a strategy consisting of $(a_0,a_1,\ldots,a_m)$, where $a_i\in U$ and are $F_i$ measurable and $(\hat{\prob}_n^{V_m})^m(x,{\o})$ stands for a controlled probability under which in the step $m$ we enter the open set ${\o}$, not leaving the set $B_n$ in the meantime.

We want to show the following theorem
\begin{theorem}\label{thm3} Under assumption \eqref{A1} and \eqref{A2'} for each $n\in \bN$ there exists $w_n\in C(E)$ and a constant $\lambda_n$ such that
\begin{equation}\label{meqBelln}
e^{w_n(x)}=\inf_{a\in U} \left[e^{c_n(x,a)-\lambda_n}\int_{E} e^{w_n(y)}\tilde{\prob}_n^a(x,dy)\right]
\end{equation}
is satisfied and $\sup_{x\in B_{n+1}}w_n(x)=0$.
\end{theorem}
\begin{proof}
We follow the arguments of the proof of Theorem \ref{thm1}. For $f\in C(E)$  we define
\begin{equation}
T_nf(x):= \inf_{a\in U}\left[e^{c_n(x,a)}\int_E f(y)\tilde{\prob}_n^a(x,dy)\right].
\end{equation}
For $x\in B_{n+1}^c$ we clearly have that $T_nf(x)=e^{\underline{c}}\mu(f)$. For such operator $T_n$ the claim of Proposition \ref{prop1} together with estimation \eqref{equicont}  is satisfied. Therefore by the Krein Rutman Theorem \ref{thmA1} (compare to Corollary \ref{cor1}) there is $v_n\in C_+(E)$ with $\|v_n\|=1$ and positive $\lambda_n$ such that for $x\in E$ and $n\in \bN$
\begin{equation}\label{vnG}
T_nv_n(x)=e^{\lambda_n} v_n(x).
\end{equation}
Using \eqref{A2'} analogously as in the proof of Lemma \ref{lem1} (see \eqref{vn positive}) we obtain that $v_n(x)>0$ for $x\in B_n$. Moreover $v_n(x)=e^{-\lambda_n}e^{\underline{c}}\mu(v_n)>0$ for $x\in B_{n+1}^c$. Let $G=\left\{x: v_n(x)=0\right\}$. Clearly $G\subset B_{n+1}\setminus B_n$.  Therefore from \eqref{vnG} we have that $\inf_{a\in U} \tilde{\prob}_n^a(x, G)=1$ for $x \in G$, but $\tilde{\prob}_n^a(x, G)=0$, a contradiction. Consequently $w_n(x)=\ln v_n(x)$ is well defined and is a continuous bounded function which together $\lambda_n$ form a solution to \eqref{meqBelln}. In addition since $\|v_n\|=1$ and $v_n$ is nonnegative we have that $\sup_{x\in B_{n+1}}w_n(x)=0$.
\end{proof}
Adapting the proof of Theorem \ref{thm3} to the case when $E$ is compact we easily obtain
\begin{proposition}\label{propex2}
In the case when $E$ is compact, under \eqref{A1} as well as the following version of \eqref{A2'}
\begin{enumerate}
\item[(\namedlabel{A2'*}{A.2'*})] for each  $x\in E$ and open set $\o$ there is $m\in \bN$  such that
    $$\inf_{V_m}({\prob}^{V_m})^m(x,{\o})> 0,$$
\end{enumerate}
there is $w\in C(E)$ and constant $\lambda$ for which equation \eqref{meqBell} is satisfied.
\end{proposition}

Our purpose in now to let $n\to \infty$ in \eqref{meqBelln} to obtain a solution to \eqref{meqBell}. Define $\check{c}_n(x)=\inf_{a\in U} c_n(x,a)$ and $\check{c}(x)=\inf_{a\in U} c(x,a)$. We shall need the following versions of the assumptions  \eqref{A3} and \eqref{A4}
\begin{enumerate}
\item[(\namedlabel{A3'}{A.3'})]  for $\lambda:=\limsup_{n \to \infty} \lambda_n$ there is $\vep>0$ such that

    the set $\left\{x: \check{c}(x)\geq \lambda-\vep\right\}$ is compact

\item[(\namedlabel{A4'}{A.4'})]  there is $x\in E$ and $m\in \bN$, $N\in \bN$ and $\vep>0$ such that $\inf_{n\geq N} \inf_{V_m}(\hat{\prob}_n^{V_m})^m(x,{\cal O}_n(\vep))>0$
 with $ {\cal O}_n(\vep)=\left\{y\in E: e^{w_n(y)}>\vep\right\}$.
\end{enumerate}
\begin{remark}
By analogy to Remark \ref{rem1} assumption \eqref{A3'} is satisfied when
 \[\limsup_{\rho(b,x)\to \infty} \check{c}(x)<\lambda,\] since then there exists $\vep>0$ and $M>0$ such that $\rho(b,x)\geq M$ implies that $\check{c}(x)<\lambda-\vep$.
By Remark \ref{rem3} under \eqref{addass} and \eqref{addass'} we have that $B_n\subset {\o}_n(e^{-K})$. Assuming additionally that there is $N\in \bN$ and $x\in E$ such that $\inf_{a\in U} \prob^a(x,B_N)>0$ we have that $\inf_{n\geq N} \inf_{V_m}(\hat{\prob}_n^{V_m})^m(x,{\cal O}_n(e^{-K}))>0$ and \eqref{A4'} is satisfied.
\end{remark}
We have
\begin{theorem}\label{thm4} Under \eqref{A1}, \eqref{A2'} when $E$ is not compact and either \eqref{A3'} or \eqref{A4'} are satisfied, there is a continuous function $w$ such that $\sup_{x\in E}w(x)=0$ and the pair $w,\lambda$ forms a solution to the equation \eqref{meqBell}.
\end{theorem}
\begin{proof}
We follow the proof of Theorem \ref{thm2}. Namely, in this case we also have \eqref{o1}-\eqref{o3} and obtain an analog of \eqref{o4}
\begin{equation}\label{o4'}
 \lim_{k\to \infty} e^{w_{n_k}(x)}= \inf_{a\in U} e^{c(x,a)-\lambda}\int_{E} z(y){\prob}^a(x,dy)=z(x).
 \end{equation}
 where $z(x)$ is a limit of $e^{w_{n_k}(x)}$ as $k\to \infty$  with convergence uniform on compact subsets and $\lim_{k\to \infty}\lambda_{n_k}=\lambda$. Notice that $w_n\not\equiv 0$ since otherwise from \eqref{meqBelln} we have that $w_n(x)=\inf_{a\in U} c_n(x,a)-\lambda_n=0$ for $x\in E$, and $\check{c}_n(x)=\lambda_n$ for $x\in E$. Consequently $\check{c}_n(x)=\underline{c}=\lambda_n$ for $x\in E$. Since $\check{c}_n(x)\leq \check{c}(x)$ for $x\in E$ we obtain $\check{c}(x)\geq \lambda - |\lambda-\lambda_n|$ for $x\in E$. Therefore whenever $|\lambda_{n_k}-\lambda|\leq \vep$ with $\vep$ from \eqref{A3'} we have a contradiction since the set $\left\{x: \check{c}(x)\geq \lambda - \vep\right\}=E$. In what follows we shall assume that $w_{n_k}\not\equiv 0$.
 Let $\hat{x}_n$ be such that $w_n(\hat{x}_n)=0$. Clearly $\hat{x}_n\in B_{n+1}$.  Then from \eqref{meqBelln} we have
 \begin{equation}\label{o5'}
 1=\inf_{a\in U} \left[e^{c_n(\hat{x}_n,a)-\lambda_n}\int_E e^{w(y)}\prob_n^a(\hat{x}_n,dy)\right]\leq e^{\check{c}_n(\hat{x}_n)-\lambda_n}
 \end{equation}
 and therefore $\check{c}_n(\hat{x}_n)\geq\lambda_n$ and $\hat{x}_n\in \left\{x: \check{c}(x) \geq \lambda_n\right\}\subset \left\{x: \check{c}(x)\geq  \lambda -(\lambda-\lambda_n)\right\}$.
 Consequently for a sufficiently large $k$, and subsequence $(n_k)$ we have that $\hat{x}_{n_k}\in \left\{x: \check{c}(x)\geq  \lambda -\vep\right\}$, with $\vep>0$ from \eqref{A3'}. This means that $\hat{x}_{n_k}$ are in a fixed compact set for $k$ sufficiently large and consequently $\|z\|=1$.
 Using assumption \eqref{A2'} in a similar way as in the proof of Theorem \ref{thm2} we obtain that $z(x)>0$ for $x\in E$, which completes the proof in the case of assumption \eqref{A3'}.
 When \eqref{A4'} is satisfied by analogy to \eqref{prA3} we have
 \begin{eqnarray}\label{prA3'}
&& v_n(x)=\inf_V \tilde{\ee}_x^V\left[\exp\left\{\sum_{i=0}^{m-1}(c_n(X_i,a_i)-\lambda_n))\right\}
v_n(X_m)\right] \geq \nonumber \\
&& \vep e^{m\underline{c} - m \lambda_n}\inf_{V_m} (\hat{\prob}_n^{V_m})^m(x,{\cal O}_n(\vep))
\end{eqnarray}
so that from \eqref{A4'} we obtain that $z(x)>0$ for $x\in B_n$, and consequently for any $x\in E$. The proof under \eqref{A4'} is therefore completed.
  \end{proof}

We have the following analog of Proposition \ref{prop2}.
\begin{proposition}\label{prop3}
Assume \eqref{A1}, \eqref{A2'}, and \eqref{A3'} or \eqref{A4'}. Then
\[\lambda\leq \inf_{x\in E}\inf_V J_x^m(V).\] If for a given $x\in E$ and any $\vep>0$ there is $n\in \bN$ such that for strategy $\hat{V}=(\hat{u}(X_i))$, where $\hat{u}$ is a Borel measurable selector of the right hand side of the equation \eqref{meqBell}, condition \eqref{exass} is satisfied,
then $\lambda=\inf_V J_x^m(V)$ and $\hat{V}$ is  an optimal control.
\end{proposition}
\begin{proof}
Iterating \eqref{meqBell} we obtain that
\begin{equation}\label{itmBell}
e^{w(x)}=\inf_V \ee_x^V\left[\exp\left\{\sum_{i=0}^{m-1}(c(X_i,a_i)-\lambda)+w(X_m)\right\}\right].
\end{equation}

Therefore  taking into account that $w_n(x)\leq 0$ for $x\in E$ for any strategy $V$ we have
\begin{eqnarray}\label{o6'}
&&\lambda\leq{1\over m} \ln \left\{\ee_x^{{V}}\left[\exp\left\{\sum_{i=0}^{m-1}c(X_i,a_i)+w(X_m)\right\}\right]\right\}
-{w(x)\over m}\leq \nonumber \\
&&{1\over m} \ln \left\{\ee_x^{{V}}\left[\exp\left\{\sum_{i=0}^{m-1}c(X_i,a_i))\right\}\right]\right\}
-{w(x)\over m}
\end{eqnarray}
and letting $m\to \infty$ we obtain that $\lambda\leq J_x^m(V)$.
On the other hand under \eqref{exass} for strategy $\hat{V}$ we have
\begin{eqnarray}\label{o7'}
&&\lambda={1\over m} \ln \left\{\ee_x^{\hat{V}}\left[\exp\left\{\sum_{i=0}^{m-1}c(X_i,\hat{u}(X_i))+w(X_m)\right\}\right]\right\}
-{w(x)\over m}\geq \nonumber \\
&&{1\over m} \ln \left\{\ee_x^{\hat{V}}\left[1_{B_n}(X_m)\exp\left\{\sum_{i=0}^{m-1}c(X_i,\hat{u}(X_i)))\right\}
e^{\inf_{y\in B_n}w(y)}+ \right. \right. \nonumber \\
&&\left. \ee_x^{\hat{V}}\left[1_{B_n^c}(X_m)e^{-m\|c\|+w(X_m)}\right]\right\}-{w(x)\over m}\geq \nonumber \\
&& {1\over m} \ln \left\{\ee_x^{\hat{V}}\left[1_{B_n}(X_m)\exp\left\{\sum_{i=0}^{m-1}c(X_i,\hat{u}(X_i)\right\}
\right]\right\}+
{\inf_{y\in B_n}w(y)\over m}-{w(x)\over m}.
\end{eqnarray}
Letting $m\to \infty$ and using \eqref{exass} we obtain
\begin{equation}
\lambda\geq \limsup_{m\to \infty}{1\over m} \ln \left\{\ee_x^{{\hat{V}}}\left[\exp\left\{\sum_{i=0}^{m-1}c(X_i,a_i))\right\}\right]\right\}-\vep.
\end{equation}
Since $\vep$ can made arbitrarily small, taking into account \eqref{o7'} we have that $\lambda=J_x(\hat{V})$, which completes the proof.

\end{proof}
\begin{remark}
All comments concerning uniqueness of solution to the equation \eqref{eqBell} pointed out in Remark \ref{remm} concern also the equation \eqref{meqBell}. Namely, when \eqref{exass} is not satisfied  $\lambda$ may not be defined in a unique way. Furthermore in such case $\lambda$ is only a lower bound for optimal values of the functional \eqref{mfun1}. The same Example 1 of \cite{Jas2007} shows that optimal values of the functional \eqref{mfun1} may depend on initial state of $(X_n)$ and are not solutions to \eqref{meqBell} for any function $w$.
\end{remark}

\section{Appendix}
We recall here Krein Rutman theorem (see original paper \cite{KreinRut} or  \cite{Bonsall}, \cite{Arap}) in an ordered Banach space $X$.
Assume that $K$ is a nontrivial closed cone with nonempty interior in $X$ such that $K\cap (-K)={0}$. This cone introduces an order $\succeq$ in $X$  such that for $x,y\in X$ we have $x\succeq y$ whenever $x-y\in K$.
\begin{theorem}\label{thmA1}
Assume $T:X\mapsto X$ is a continuous map such that $T(K)\subset K$, which is order preserving i.e. whenever $x\succeq y$ we have $Tx\succeq Ty$, and $1$-homogeneous, i.e. $T(a x)=aTx$ for $x\in X$ and $a\geq 0$. Furthermore assume that $T$ is completely continuous in $X$ and for some $y\in K$ there is $M>0$ such that $MTy\succeq y$. Then there exist $\lambda>0$ and $\bar{x}\in K$ such that $\|\bar{x}\|=1$ and $T\bar{x}=\lambda \bar{x}$.
\end{theorem}

\bibliographystyle{siamplain}

\end{document}